\documentclass[11pt]{article}
\usepackage{amsmath,amstext, verbatim,remreset}
\usepackage[T2A]{fontenc}
\usepackage[utf8]{inputenc}
\usepackage[ukrainian]{babel}
\usepackage{amsfonts,amssymb}
\usepackage{cite}

\newcommand{\RR}{\mathbb R}
\newcommand{\NN}{\mathbb N}

\newcommand{\K}{{\bf k}}
\usepackage{mathrsfs}
\newtheorem{theorem}{Теорема}
\newtheorem{lemma}{Лема}
\newtheorem{definition}{Означення}
\newtheorem{problem}{Задача}
\newtheorem{remark}{Зауваження}
\begin{document}
\begin{center}
{\bf \Large Про задачу Колмогорова на класах абсолютно монотонних функцій, кратно монотонних функцій, та її зв'язок з проблемою моментів Маркова}

\bigskip
{\bf \Large В. Ф. Бабенко, Ю. В. Бабенко, О. В. Коваленко}
\end{center}
\section{Постановка задачі}
Для $r\in\NN$ покладемо $L^r_{\infty,\infty}(\RR_-):=L_\infty(\RR_-)\bigcap L^r_{\infty}(\RR_-),$ $\|\cdot\|:=\|\cdot\|_{L_\infty(\RR_-)}$.
Ми розглядаємо задачу Колмогорова у наступній постановці.

\begin{problem}[Задача Колмогорова]
	Нехай задано клас функцій $X\subset L^r_{\infty,\infty}(\RR_-)$ і довільна система $d$ цілих чисел
	$0\leq k_1 < k_2 < \ldots < k_{d}\leq r.$
	Знайти необхідні та достатні умови на систему додатних чисел $M_{k_1}, M_{k_2}, \dots, M_{k_{d}}$
	для того, щоб гарантувати існування функції $x \in X$ такої, що $$\|x^{(k_i)}\|=M_{k_i},\; i=1,\dots,d.$$
\end{problem}

Для $d\in\NN$ і цілих чисел $0\leq k_1<k_2<\ldots<k_{d}\leq r$ покладемо ${\bf k}:=(k_1,\ldots, k_{d})$, $M_{{\bf k}} := \{M_{k_1},\ldots ,M_{k_{d}}\}.$

Для функції $x\in X$ покладемо $$M_{\bf k}(x):=(M_{k_1}(x),\ldots ,M_{k_d}(x)),$$ де $$M_{k_i}(x)=\| x^{(k_i)}\|,\;\; i=1,\ldots ,d.$$

\begin{definition}
Будемо називати набір чисел $M_{{\bf k}}$ {\bf допустимим} для класу $X\subset L^r_{\infty,\infty}(\RR_-)$, якщо існує функція $x\in X$  така, що $\left\|x^{(k_i)}\right\|=M_{k_i}$, $i=1,2,\dots,d$
(або більш коротко $M_{{\bf k}}(x)=M_{{\bf k}}$). 
\end{definition}

Множину всіх ненульових допустимих наборів $M_{{\bf k}}$ ми будемо позначати через $A_{\K}(X)$.

В наведених позначеннях задачу Колмогорова можна сформулювати наступним чином.

Для заданого класу функцій $X\subset L^r_{\infty,\infty}(\RR_-)$ і довільної системи $d$ цілих чисел $\K$ охарактеризувати множину $A_{\K}(X)$.

Ми також розглядаємо задачу Колмогорова в дещо іншому формулюванні.
\begin{problem}[Задача Колмогорова (альтернативне формулювання)]
Для класу функцій $$X\subset L^r_{\infty,\infty}(\RR_-)$$ і довільної системи $d$ цілих чисел $\K$ знайти ''мінімальну'' множину $F_{\K}(X)\subset X$ таку, що $$A_{\K}(X) = \left\{M_{\K}(x)\colon x\in F_{\K}(X)\right\}.$$
\end{problem}

Історію питання та огляд відомих результатів можна знайти в статті~\cite{Kovalenko13a}.

\section{Допоміжні результати}

В цьому пункті ми наведемо необхідні нам у подальшому означення та результати, пов'язані з розв'язком проблеми моментів Маркова. Класичні результати у цій проблематиці можна знайти у монографіях~\cite{Karlin66,Akhiezer65,Krein73}. 

Означення та результати, які містяться в даному пункті, сформульовані у термінах і формулювання близьких до тих, що використовуються у монографії~\cite{Karlin66}, глави~2 та~5.

\begin{definition}\label{tchebishovSystem}
	Система функцій $u_1,\dots, u_n$ називається чебишовською на $[0,\infty)$, якщо функції $u_1,\dots, u_n$ неперервні на $[0,\infty)$ і 
	$$\det \left[\|u_i(t_j)\|_{i,j=1}^n\right]>0$$ при довільному виборі точок $0\leq t_1<\ldots<t_n<\infty$.
\end{definition}

\begin{problem}[Проблема моментів Маркова]\label{markovProblem} 
	Нехай задано чебишовську на $[0,\infty)$ систему функцій $u_1,\dots,u_n$. Знайти необхідні та достатні умови на набір чисел 
	${\bf c} = (c_1,\dots,c_n)\in\RR^{n}$, для того, щоб гарантувати існування функції $\sigma$ з множини $\mathcal{D}$ невід'ємних 
	неспадних функцій обмеженої варіації, такої, що 
	\begin{equation}\label{moments}
		c_k = \int\limits_0^\infty u_k(t)d\sigma(t),\,k=1,\dots,n.
	\end{equation}
\end{problem}

\begin{remark} 
	У формулюванні задачі вище розглядаються тільки ті функції $\sigma\in\mathcal{D}$, для яких інтеграли в~\eqref{moments} є 
	абсолютно збіжними.
\end{remark}

\begin{definition}\label{momentSpace} 
	Позначимо через $$\mathcal{M}_n = \mathcal{M}(u_1,\dots,u_n)$$ множину всіх точок $(c_1,\dots,c_n)\in\RR^{n}$,
	для яких існує функція $\sigma\in\mathcal{D}$ така, що виконуються рівності~\eqref{moments}.
\end{definition}
Справедлива наступна теорема. 

\begin{theorem}\label{conicalHull}
	Нехай на $[0,\infty)$ задано систему неперервних функцій $\left\{u_k\right\}_{k=1}^n$ (не обов'язково чебишовську). 
	Тоді моментний простір $\mathcal{M}(u_1,\dots,u_n)$ співпадає з опуклим конусом, натягнутим на підмножину 
	$$\mathcal{C}_n = \mathcal{C}(u_1,\dots,u_n):=\left\{(u_1(t),\dots,u_n(t)),\, t\in [0,\infty)\right\}$$ 
	моментного простору $\mathcal{M}(u_1,\dots,u_n)$, що породжується точковими мірами.
\end{theorem}

З теореми~\ref{conicalHull} і теореми Каратеодорі слідує, що кожна точка з множини $\mathcal{M}(u_1,\dots,u_n)$ може бути представлена у вигляді лінійної комбінації не більше ніж $n+1$ точок кривої $\mathcal{C}(u_1,\dots,u_n)$.

Таким чином, довільну точку ${\bf c}\in \mathcal{M}_n$ можна представити у вигляді
\begin{equation}\label{representation}
	{\bf c} = \sum\limits_{k=1}^ma_ku(t_k),\;m\leq n+1,
\end{equation}
де $0\leq t_1<\ldots<t_m$, $a_k>0$, $k=1,\dots,m$ і $u(t) := (u_1(t), u_2(t),\dots,u_n(t))$.

\begin{definition}
	Числа $t_k$, $k=1,\dots,m$ з рівності~\eqref{representation} ми будемо називати коренями представлення~\eqref{representation}.
\end{definition}

\begin{remark}\label{uConditions}
	У випадку, коли система функцій $\left\{u_k\right\}_{k=1}^n$  задовольняє деяким додатковим умовам, представлення~\eqref{representation} може містити і менше ніж $n+1$ точку. Достатніми для цього умовами є наступні:
	\begin{enumerate}
		\item\label{uCondition1} Системи функцій $\left\{u_k\right\}_{k=1}^n$ і $\left\{u_k\right\}_{k=1}^{n-1}$ є чебишовськими на $[0,\infty)$.
		\item Існує поліном $u(t)=\sum\limits_{k=1}^na_ku_k(t)$, що задовольняє умовам $u(t)>0$, $t\in [0,\infty)$ і 
			$\varliminf\limits_{t\to\infty} u(t) > 0$.
		\item $\lim\limits_{t\to\infty}\frac{u_k(t)}{u_n(t)} = 0$, $k=1,\dots,n-1$.
	\end{enumerate}
\end{remark}

Нам знадобиться наступне означення.

\begin{definition}\label{index} 
	Індексом $I({\bf c})$ точки ${\bf c}\in \mathcal{M}_n$ ми будемо називати мінімальну кількість точок з $\mathcal{C}_n$, які можуть представити 
	${\bf c}$ як опуклу комбінацію. При цьому точки $(u_1(t),\dots,u_n(t))$ з $t > 0$ рахуються за одиницю, а з $t = 0$ --- за половину.
\end{definition}

Якщо функції $u_1,\dots, u_n$ задовольняють умовам з зауваження~\ref{uConditions}, то справедливі наступні теореми.
\begin{theorem}\label{momentSpaceBoundary}
	Ненульовий вектор ${\bf c}$ належить границі $\partial\mathcal{M}_n \bigcap \mathcal{M}_n$ моментного простору $\mathcal{M}_n$ тоді і 
	тільки тоді, коли $I({\bf c}) < \frac n 2$.
\end{theorem}

\begin{theorem}\label{momentSpaceMainRepresentation}
	Нехай ${\bf c}\in {\rm int}\,\mathcal{M}_n$. Тоді існує представлення ${\bf c}$ індексу $\frac n 2$.
\end{theorem}

\begin{theorem}\label{momentSpaceCanonicalRepresentation}
	Нехай ${\bf c}\in {\rm int}\,\mathcal{M}_n$. Тоді для всіх $t^*>0$ існує представлення вектора ${\bf c}$ індексу $\frac {n + 1} 2$,
	коренем якого є число $t^*$.
\end{theorem}

\begin{remark}
	Нехай задано $d\in\NN$ і цілі числа $0\leq k_1<k_2<\ldots<k_d$. У подальшому ми будемо розглядати тільки системи функцій 
	$\left\{u_i\right\}_{i=1}^d$ наступного вигляду: $u_i(t) = t^{k_i}$, $i=1,\dots, d$. У випадку $k_1=0$ такі системи задовольняють умовам з
	зауваження~\ref{uConditions}, а отже для них справедливі теореми~\ref{momentSpaceBoundary}, \ref{momentSpaceMainRepresentation} і~\ref{momentSpaceCanonicalRepresentation}.
\end{remark}
	Доведення справедливості умови~\ref{uCondition1} зауваження~\ref{uConditions} для системи $\left\{t^{k_i}\right\}_{i=1}^d$ у випадку $k_1=0$ міститься у~\cite[гл.~1 \S 3]{Karlin66}, дві інші умови виконуються очевидним чином.

\section{Клас абсолютно монотонних на $\RR_-$ функцій}
\begin{definition} 
	Нескінченно диференційовну на $\RR_-$ функцію будемо називати абсолютно монотонною, якщо вона і всі її похідні невід'ємні на $\RR_-$. 
	Через $AM(\RR_-)$ ми будемо позначати клас абсолютно монотонних на $\RR_-$ функцій.
\end{definition}

Справедливе наступне інтегральне представлення абсолютно монотонних функцій, що було доведено Берштейном~\cite{Bernshtein}. 

\begin{theorem}\label{AMrepresentation}
	Функція $x(t)$ є абсолютно монотонною тоді і тільки тоді, коли її можна представити за допомогою інтегрального представлення
	\begin{equation}\label{Bern}
		x(t)= \int_0^{\infty} e^{t u} d\beta(u),\; t\in\RR_-,
	\end{equation}
	де $\beta(u)$ --- неспадна обмежена функція.
\end{theorem}

Відмітимо, що, в силу означення, рівномірні норми абсолютно монотонної на $\RR_-$ функції $x(t)$ та її похідних досягаються у точці нуль. В силу інтегрального представлення~\eqref{Bern} це означає, що для $k=0,1,\dots$ 
$$\|x^{(k)}\| = x^{(k)}(0) =  \int_0^{\infty} u^k d\beta(u).$$

Таким чином, справедлива наступна теорема, яка показує зв'язок розв'язку задачі Колмогорова на класі абсолютно монотонних функцій і проблеми моментів Маркова.
\begin{theorem}\label{AMandMomentProblem}
	Нехай задано $d\in\NN$, $\K = (0\leq k_1 < k_2<\ldots<k_d)$. Тоді множина $A_{\K}(AM(\RR_-))$ допустимих для класу $AM(\RR_-)$ наборів співпадає з моментним простором $\mathcal{M}(t^{k_1},\dots,t^{k_d})$.
\end{theorem}

Відмітимо, що точкова міра $\beta$, яка  зосереджена у точці $a > 0$ породжує абсолютно монотонну функцію $e^{at}$. Нам буде зручно розглядати абсолютно монотонні функції $e^{\frac t a}$, що породжується точковими мірами в точках $a^{-1}$, $a>0$. 

\begin{definition}
	Функцію $$\phi(AM(\RR_-),{\bf a}, {\bf \lambda};t):=\sum\limits_{s=1}^m \lambda_s a_s^r e^{a_s^{-1} t},$$ де 
	$\lambda_s, a_s > 0$, $s=1,\dots,m$, ${\bf \lambda} = (\lambda_1,\dots,\lambda_m)$, ${\bf a} = (a_1,\dots,a_m)$ 
	ми будемо називати $AM(\RR_-)$ -- ідеальним сплайном порядку $r\in\NN$ з $m$ вузлами $-a_1,\dots,-a_m$.
	
	Якщо $\phi(t)$ --- $AM(\RR_-)$ -- ідеальний сплайн з $m$ вузлами і $C > 0$, то функцію $C+\phi(t)$ ми будемо називати 
	$AM(\RR_-)$ -- ідеальним сплайном з $m+\frac 1 2$ вузлами.
\end{definition}

\section{Клас кратно монотонним за Вільямсоном функцій}

Через $L_{\infty,\infty}^{r,\cup}(\RR_-)$ позначимо клас функцій $x\in L_{\infty,\infty}^{r}(\RR_-)$ таких, що для $k=0,\dots, r-1$ похідні $x^{(k)}$ є неспадними та опуклими (див.~\cite{Williamson}).


Вільямсон~\cite{Williamson} довів наступну теорему. 
\begin{theorem}
	$y(t) \in L_{\infty,\infty}^{r,\cup}(\RR_-)$ тоді і тільки тоді, коли
	\begin{equation}\label{Will}
		y(t)=\frac{1}{r!}\displaystyle \int_0^{\infty} \left[(1+ut)_+\right]^rd\beta(u),\;t\in\RR_-,
	\end{equation}
	де $\beta(u)$ --- неспадна обмежена функція.
\end{theorem}

Відмітимо, що точкова міра $\beta(t)$, що зосереджена у точці $a$, $a>0$ дає функцію  $\frac 1{r!}(1+at)^r_+$. Нам буде зручно розглядати функцію $\frac 1{a^rr!}(a+t)^r_+$, що породжується точковою мірою у точці $\frac 1 a$.

\begin{definition}
	Функцію $$\phi(L_{\infty,\infty}^{r,\cup}(\RR_-),{\bf a}, {\bf \lambda};t):=\frac 1 {r!}\sum\limits_{s=1}^m \lambda_s (a_s + t)^r_+,$$ де 
	$\lambda_s, a_s > 0$, $s=1,\dots,m$, ${\bf \lambda} = (\lambda_1,\dots,\lambda_m)$, ${\bf a} = (a_1,\dots,a_m)$ ми будемо 
	називати $L_{\infty,\infty}^{r,\cup}(\RR_-)$ -- ідеальним сплайном порядку $r\in\NN$ з $m$ вузлами $-a_1,\dots,-a_m$.
	
	Якщо $\phi(t)$ --- $L_{\infty,\infty}^{r,\cup}(\RR_-)$ -- ідеальний сплайн з $m$ вузлами і $C > 0$, то функцію $C+\phi(t)$ ми будемо називати 
	$L_{\infty,\infty}^{r,\cup}(\RR_-)$ -- ідеальним сплайном з $m+\frac 1 2$ вузлами.
\end{definition}

\section{Зв'язок між класами $L_{\infty,\infty}^{r,\cup}(\RR_-)$ і $AM(\RR_-)$}

Для чисел $a_1,\dots, a_d$ через ${\rm diag}(a_1,\dots,a_d)$ позначимо квадратну діагональну матрицю порядку $d$ з числами $a_1,\dots, a_d$ на головній діагоналі. Для даного вектора ${\bf c} \in \RR^d$ позначимо через ${\rm diag}(a_1,\dots, a_d)\,{\bf c}$ результат множення матриці ${\rm diag}(a_1,\dots,a_d)$ на вектор-стовпець ${\bf c}$. Для множини $A\subset \RR^d$ покладемо $${\rm diag}(a_1,\dots, a_d)\,A:=\left\{{\rm diag}(a_1,\dots, a_d)\,{\bf c}\colon {\bf c} \in A\right\}.$$

Враховуючи представлення \eqref{Bern}, \eqref{Will} і той факт, що усі норми похідних порядків $k_1,\dots,k_d$ у функцій з класів $L_{\infty,\infty}^{r,\cup}(\RR_-)$ і $AM(\RR_-)$ досягаються в точці нуль, ми отримуємо наступну теорему.

\begin{theorem}\label{connection}
	Нехай задано цілі числа $0\leq k_1<k_2<\ldots<k_d\leq r$, $\K = (k_1,\dots,k_d)$. Тоді 
	$$A_{\K}(AM(\RR_-)) = {\rm diag}((r-k_1)!,\dots, (r-k_d)!) A_{\K}(L_{\infty,\infty}^{r,\cup}(\RR_-)).$$
\end{theorem}

\begin{remark}\label{corresp}
	Нехай задано вектори ${\bf a}, {\bf \lambda}$ додатних чисел. Тоді 
	$$M_{\K}(\phi(AM(\RR_-),{\bf a},{\bf \lambda};t)) = {\rm diag}((r-k_1)!,\dots, (r-k_d)!) M_{\K}(\phi(L_{\infty,\infty}^{r,\cup}(\RR_-),{\bf a},{\bf \lambda};t)).$$
\end{remark}

\section{Деякі властивості ідеальних сплайнів}

Наступна лема встановлює деякі екстремальні властивості введених сплайнів.

\begin{lemma}\label{l::1}
	Нехай $x\in L_{\infty,\infty}^{r,\cup}(\RR_-)$, $0 \leq k_0 < k_1<\ldots<k_{d}=r$. Нехай $L_{\infty,\infty}^{r,\cup}(\RR_-)$ -- ідеальний сплайн $\phi(t)$ з не більш ніж $\left[\frac d 2\right]$ вузлами такий, що 
	$\left\|\phi^{(k_i)}\right\| = \left\|x^{(k_i)}\right\|$, $i = 1,\dots,2\cdot\left[\frac d 2\right]$. Тоді
	\begin{equation}\label{l1}
		\left\|\phi^{(k_0)}\right\|\leq \left\|x^{(k_0)}\right\|.
	\end{equation}
	У випадку непарного $d$ має місце нерівність
	\begin{equation}\label{l2}
		\left\|\phi^{(r)}\right\|\leq \left\|x^{(r)}\right\|.
	\end{equation}	
	Якщо в нерівності~\eqref{l1} має місце рівність, $d$ --- непарне і в~\eqref{l2} має місце рівність, або сплайн $\phi$ має менше ніж $\left[\frac d 2\right]$ вузлів, то $x^{(k_1)} \equiv \phi^{(k_1)}$.
\end{lemma}
{\bf Доведення.} Нехай $d$ парне. Тоді $2\cdot\left[\frac d 2\right] = d$. Припустимо супротивне, нехай $x^{(k_0)} \ne \phi^{(k_0)}$ і $\|x^{(k_0)}\| \leq \|\phi^{(k_0)}\|$.
Покладемо $\Delta(t):=x(t)-\phi(t)$. Для того, щоб отримати суперечність, ми будемо рахувати число змін знаку різниці $\Delta(t)$ та її похідних.

Перш за все відмітимо, що в силу означення сплайнів $\phi$, маємо $\phi^{(k_0)}(-a_1)=0$ (де $-a_1$ --- самий лівий вузол сплайна $\phi$). Крім того, $x^{(k_0)}(-a_1)\geq 0$, і отже ми отримуємо, що $\Delta^{(k_0)}(-a_1)\geq 0$. В силу припущення
$$\Delta^{(k_0)}(0)=x^{(k_0)}(0)-\phi^{(k_0)}(0)=\left\|x^{(k_0)}\right\|-\left\|\phi^{(k_0)}\right\|\leq 0.$$ 
Тоді існує точка $t_{k_0+1}^1\in(-a_1,0)$ така, що $\Delta^{(k_0+1)}(t_{k_0+1}^1)< 0$. Крім того, 
$\Delta^{(k_0+1)}(-a_1)\geq 0$. Тому існує точка $t_{k_0+2}^1\in(-a_1,0)$ така, що $\Delta^{(k_0+2)}(t_{k_0+2}^1)<0$. Повторюючи аналогічні міркування ми отримаємо, що існує точка $t_{k_1}^1\in(-a_1,0)$ така, що $\Delta^{(k_1)}(t_{k_1}^1)< 0$. Крім того, $\Delta^{(k_1)}(-a_1)\geq
0$ і $\Delta^{(k_1)}(0)= 0$ за умовою леми. Таким чином, існують точки $-a_1<t^1_{k_1+1}<t^2_{k_1+1}<0$ такі, що $\Delta^{(k_1+1)}(t^1_{k_1+1})<0$ і $\Delta^{(k_1+1)}(t^2_{k_1+1})>0$. Такий розподіл знаків збережеться до рівня $k_2$ де, враховуючи умову $\left\|x^{(k_2)}\right\|=\left\|\phi^{(k_2)}\right\|$ і той факт, що $\Delta^{(k_2)}(-a_1)\geq 0$, ми отримаємо існування точок $-a_1<t^1_{k_2+1}<t^2_{k_2+1}<t^3_{k_2+1}<0$ таких, що $\Delta^{(k_2+1)}(t^1_{k_2+1})<0$, $\Delta^{(k_2+1)}(t^2_{k_2+1})>0$, і $\Delta^{(k_2+1)}(t^3_{k_2+1})<0$. Продовжуючи аналогічно, ми отримаємо, що існують точки $-a_1<t_{k_{d-1}+1}^1<\ldots<t_{k_{d-1}+1}^{d}<0$ такі, що $(-1)^i\Delta^{(k_{d-1}+1)}(t_{k_{d-1}+1}^{i})>0$, $i=1,\dots,d$, і так далі, до рівня $r-1$. 

На рівні $r-1$-ї похідної існують точки $-a_1<t_{r-1}^1<\ldots<t_{r-1}^{d}<0$ такі, що $(-1)^i\Delta^{(r-1)}(t_{r-1}^{i})>0$, $i=1,2,\dots,d$. Крім того, $\Delta^{(r-1)}(-a_1)\geq 0$. 

Тоді на інтервалі $(-a_1, t^1_{r-1})$ існує множина $S_0$ додатної міри така, що $\Delta^{(r)}(t)<0$ для всіх $t\in S_0$. Крім того, для всіх $i=1,\dots, d-1$ на інтервалі $(t^i_{r-1}, t^{i+1}_{r-1})$ існують множини $S_i\subseteq (t^i_{r-1}, t^{i+1}_{r-1})$ додатної міри такі, що $(-1)^i\Delta^{(r)}(t)<0$ для всіх $t \in S_i$. Таким чином функція $\Delta^{(r)}(t)$ має не менше ніж $d-1$ істотних змін знаку на $(-a_1, 0)$. Але це неможливо, оскільки сплайн $\phi$ має не більше ніж $\frac d 2$ вузлів і функція $\Delta^{(r)}(t)$ може змінювати знак в вузлах сплайна $\phi$ і не більше одного разу на кожному з інтервалів між вузлами сплайна $\phi$ (і, крім того, за умовою леми  $\left\|\phi^{(r)}\right\|= \left\|x^{(r)}\right\|$). Ми отримали суперечність, а отже лему у випадку парного $d$ доведено.

Нехай тепер $d$ непарне. Тоді $2\cdot\left[\frac d 2\right] = d - 1$. Проводячи міркування аналогічно до наведених вище, ми отримаємо справедливість~\eqref{l2}, причому рівність нерівності~\eqref{l2} можлива лише при $x^{(k_1)} \equiv \phi^{(k_1)}$. Якщо ж $\left\|\phi^{(r)}\right\|< \left\|x^{(r)}\right\|$, то знову  проводячи міркування аналогічно до наведених вище, ми отримаємо справедливість нерівності~\eqref{l1}, причому зі строгим знаком нерівності. {\bf Лему доведено}.

Враховуючи теореми~\ref{AMandMomentProblem}, \ref{connection} і зауваження~\ref{corresp} з теорем~\ref{momentSpaceBoundary}, \ref{momentSpaceMainRepresentation} і \ref{momentSpaceCanonicalRepresentation} ми отримуємо наступні леми.

\begin{lemma}\label{boundary}
	Нехай задано цілі числа $0 = k_1<k_2<\ldots<k_d\leq r$, $\K = (k_1,\dots,k_d)$. Нехай  $X$ позначає один з класів  
	$L_{\infty,\infty}^{r,\cup}(\RR_-)$ або $AM(\RR_-)$. ${M_{\K}}\in\partial A_{\K}(X)\bigcap A_{\K}(X)$  тоді і тільки тоді, 
	коли існує $X$ -- ідеальний сплайн $\phi(t)$ з не більше ніж $\frac {d-1} 2$ вузлами такий, що 
	\begin{equation}\label{*1}
		M_{\K}(\phi) = M_{\K}.
	\end{equation}
\end{lemma}

\begin{lemma}\label{mainRepresentation}
	Нехай задано цілі числа $0 = k_1<k_2<\ldots<k_d\leq r$, $\K = (k_1,\dots,k_d)$ і $X$ позначає один з класів  
	$L_{\infty,\infty}^{r,\cup}(\RR_-)$ або $AM(\RR_-)$. Нехай також ${M_{\K}}\in {\rm int}\,A_{\K}(X)$. Тоді існує 
	$X$ -- ідеальний сплайн $\phi(t)$ з  $\frac d 2$ вузлами такий, що виконується~\eqref{*1}.
\end{lemma}

\begin{lemma}\label{canonicalRepresentation}
	Нехай задано цілі числа $0 = k_1<k_2<\ldots<k_d\leq r$, $\K = (k_1,\dots,k_d)$ і $X$ позначає один з класів  
	$L_{\infty,\infty}^{r,\cup}(\RR_-)$ або $AM(\RR_-)$. Нехай також ${M_{\K}}\in {\rm int}\,A_{\K}(X)$. Тоді для всіх $a^*>0$ існує 
	$X$ -- ідеальний сплайн $\phi(t)$ з  $\frac {d+1} 2$ вузлами, один з яких знаходиться в точці $a^*$, такий, що виконується~\eqref{*1}.
\end{lemma}

Справедливі наступні леми.
\begin{lemma}\label{l::0.0}
	Нехай задано цілі числа $0<k_1<k_2<\ldots<k_{d}\leq r$, $\K = (k_1,\dots,k_{d})$ і $X$ позначає один з класів  $L_{\infty,\infty}^{r,\cup}(\RR_-)$ або $AM(\RR_-)$. Нехай також 
	$M_{\K}\in \partial A_{\K}(X)\bigcap A_{\K}(X)$. 
	Тоді існує $X$ -- ідеальних сплайн $\phi(t)$ з $s\leq\left[\frac{d - 1}{2}\right]$ вузлами, такий, що виконується~\eqref{*1}.
\end{lemma}

{\bf Доведення.} Нехай функція $x\in X$ така, що $M_{\K}(x) = M_{\K}$. 

Нехай $d=2n+1$. Згідно з лемою~\ref{boundary} існує $X$ -- ідеальний сплайн $\phi(t)$ з $s\leq n = \left[\frac{d-1}{2}\right]$ і $C \geq 0$ такі, що 
$M_{\bf K}(\phi+C) = M_{\bf K}(x)$, де ${\bf K} = (0,k_1,\dots, k_d)$. Це означає, що $M_{\K}(\phi) = M_{\K}(x)$ і сплайн $\phi(t)$ є шуканим.

Нехай тепер $d = 2n$. Тоді згідно з лемою~\ref{boundary} існує $X$ -- ідеальний сплайн $\phi(t)$ з $s\leq n$ вузлами такий, що $M_{\bf K}(\phi) = M_{\bf K}(x)$. Якщо $s = n$, то для всіх $C > 0$, згідно з лемою~\ref{boundary}, $M_{\bf K}(\phi+C)\in {\rm int} A_{\bf K}(X)$, що суперечить умові $M_{\K}\in \partial A_{\K}(X)\bigcap A_{\K}(X)$. Таким чином $s\leq n-1 = \left[\frac{d-1}{2}\right]$. Лему доведено.

\begin{lemma}\label{l::0}
	Нехай задано цілі числа $0\leq k_1<k_2<\ldots<k_{2d}\leq r$, $\K = (k_1,\dots,k_{2d})$ і $X$ позначає один з класів  $L_{\infty,\infty}^{r,\cup}(\RR_-)$ або $AM(\RR_-)$. Нехай також $M_{\K}\in {\rm int} A_{\K}(X)$. 
	Тоді існує $X$ -- ідеальних сплайн $\phi(t)$ з $d$ вузлами, такий, що виконується~\eqref{*1}.
\end{lemma}
{\bf Доведення.} У випадку, коли $k_1 = 0$, твердження леми одразу слідує з леми~\ref{mainRepresentation}.

Нехай тепер $k_1>0$. Розширимо систему $\{t^{k_i}\}_{i=1}^{2d}$ до системи $\{t^{k_i}\}_{i=0}^{2d}$, де $k_0 = 0$. Нехай для деякої функції $x(t)\in X$ має місце $M_{\K}(x)= M_{\K}\in {\rm int} A_{\K}(X)$. Тоді для всіх $\varepsilon > 0$ $$\left(\left\|x\right\|+\varepsilon,\left\|x^{(k_1)}\right\|,\dots,\left\|x^{(k_{2d})}\right\|\right)\in {\rm int} A_{\bf K}(X),$$ де ${\bf K} = (0,k_1,k_2,\dots,k_{2d})$. Застосовуючи лему~\ref{mainRepresentation}, отримаємо існування чисел $C\geq 0$ і $X$ -- ідеального сплайну $\phi(t)$ з $d$ вузлами такого, що  $M_{\bf K}(\phi + C) = M_{\bf K}(x + \varepsilon)$. Це означає, що виконується рівність~\eqref{*1}. Теорему доведено.

\begin{lemma}\label{l::2}
	Нехай задано цілі числа $0\leq k_1<k_2<\ldots<k_{2d+1}\leq r$, $\K = (k_1,\dots,k_{2d+1})$ і $X$ позначає один з класів  $L_{\infty,\infty}^{r,\cup}(\RR_-)$ або $AM(\RR_-)$. Нехай також $M_{\K}\in {\rm int} A_{\K}(X)$. 
	Тоді для всіх $a^*>0$ існує $X$ -- ідеальний сплайн $\phi(t)$ з  $d+1$ вузлом, один з яких знаходиться в точці $a^*$ такий, що виконується~\eqref{*1}.
\end{lemma}
 Доведення цієї леми аналогічне до доведення леми~\ref{l::0}, тільки замість леми~\ref{mainRepresentation} слід застосовувати лему~\ref{canonicalRepresentation}.

\begin{remark}\label{z::2}
	Нехай $X$ позначає один з класів $L_{\infty,\infty}^{r,\cup}(\RR_-)$ або $AM(\RR_-)$.  В умовах леми~\ref{l::0} $($або леми~\ref{l::0.0}$)$ через $\phi(X,M_{\K};t)$ будемо позначати 
	$X$ -- ідеальний сплайн з $d$ $\big($відповідно з не більш ніж $\left[\frac{d - 1}{2}\right]$ $\big)$ вузлами, для якого виконується~\eqref{*1}.
\end{remark}

\begin{remark}\label{z::3}
	З леми~\ref{l::1} слідує, що сплайн $\phi(X,M_{\K};t)$ (визначений у зауваженні~\ref{z::2}) є єдиним.
\end{remark}

Для доведення основної теореми нам знадобиться наступна проста лема.

\begin{lemma}\label{l::-1}
	Нехай задано числа $\alpha>\beta>0$, $\varepsilon>0$ і функція $\lambda(t)\colon\RR_+\to\RR_+$. Якщо для всіх достатньо великих $t>0$ $\lambda(t)\cdot t^\beta > \varepsilon$, 
	то $\lim\limits_{t\to+\infty}\lambda(t)\cdot t^\alpha=+\infty$.
\end{lemma}

\section{Розв'язок задачі Колмогорова на класах $AM(\RR_-)$ та $L_{\infty,\infty}^{r,\cup}(\RR_-)$}

З теореми~\ref{conicalHull}, лем~\ref{boundary} і~\ref{mainRepresentation} та зв'язку між класами  $AM(\RR_-)$ і $L_{\infty,\infty}^{r,\cup}(\RR_-)$ (теорема~\ref{connection}) одразу отримуємо розв'язок задачі Колмогорова для класів $AM(\RR_-)$ і $L_{\infty,\infty}^{r,\cup}(\RR_-)$ у альтернативній формі.

\begin{theorem}\label{alternativeSolution}
	Нехай задано $d\in\NN$ і цілі числа $0= k_1<k_2<\ldots<k_{d}\leq r$, $\K = (k_1,\dots,k_{d})$. Нехай $X$ позначає один з класів $L_{\infty,\infty}^{r,\cup}(\RR_-)$ або $AM(\RR_-)$. Тоді 
	$$F_{\K}(X) = \bigg\{C+\phi(X,{\bf a},{\bf \lambda}),{\bf \lambda} = (\lambda_1,\dots, \lambda_{m}), \lambda_1,\dots, \lambda_{m} > 0,$$
	$${\bf a} = (a_1,\dots, a_{m}), a_1>a_2>\ldots>a_{m}>0, C\geq 0, m = \left[\frac d 2\right]\bigg\}.$$
	У випадку парного $d$ параметр $C$ може бути покладеним рівним нулю.
\end{theorem}

Для цілочисельного вектора ${\bf k} = (k_1,\dots,k_d)$, $0\leq k_1<k_2<\ldots<k_d\leq r$ ми будемо використовувати наступне позначення: $^2{\bf k}^2 = (k_2,\dots,k_{d-1})$ (при цьому, як і раніше, ${\bf k}^2 = (k_2,\dots,k_{d})$).

\begin{theorem}\label{main}
	Нехай $d\in\NN$, $d\geq 3$ і $0\leq k_1<k_2<\ldots<k_{d}=r$ --- невід'ємні цілі числа, ${\bf k} = (k_1,\dots,k_d)$. Нехай $X$ позначає один з класів $L_{\infty,\infty}^{r,\cup}(\RR_-)$ або $AM(\RR_-)$. У випадку, коли $d$ непарне,
	$$\{ M_{{\bf k}} \in A_{{\bf k}}(X)\}\Longleftrightarrow
		\left\{\begin{array}{c}
			M_{{\bf k}^2} \in {\rm int} A_{\K^2}(X)\\
			M_{k_1}\ge \|\phi^{(k_1)}(X,M_{{\bf k}^2})\|\\
		\end{array}\right\}
		\bigvee 
	$$
	$$
		\bigvee 
		\left\{\begin{array}{c}
			M_{{\bf k}^2} \in \partial A_{\K^2}(X)\bigcap A_{\K^2}(X)\\
			k_1>0\\
			M_{k_1}=\| \phi^{(k_1)}(X,M_{{\bf k}^2})\| \\
		\end{array}\right\}
		\bigvee
		\left\{\begin{array}{c}
			M_{{\bf k}^2} \in \partial A_{\K^2}(X)\bigcap A_{\K^2}(X)\\
			k_1=0\\
			M_{k_1}\ge \|\phi^{(k_1)}(X,M_{{\bf k}^2})\|\\
		\end{array}\right\},
	$$
	а у випадку, коли $d$ парне,
	$$\{ M_{{\bf k}} \in A_{{\bf k}}(X)\}\Longleftrightarrow
		\left\{\begin{array}{c}
			M_{{\bf k}^2} \in {\rm int} A_{\K^2}(X)\\
			M_{k_1} > \|\phi^{(k_1)}(X,M_{^2{\bf k}^2})\|\\
		\end{array}\right\}
		\bigvee 
	$$
	$$
		\bigvee 
		\left\{\begin{array}{c}
			M_{{\bf k}^2} \in \partial A_{\K^2}(X)\bigcap A_{\K^2}(X)\\
			k_1>0\\
			M_{k_1}=\| \phi^{(k_1)}(X,M_{^2{\bf k}^2})\| \\
		\end{array}\right\}
		\bigvee
		\left\{\begin{array}{c}
			M_{{\bf k}^2} \in \partial A_{\K^2}(X)\bigcap A_{\K^2}(X)\\
			k_1=0\\
			M_{k_1}\ge \|\phi^{(k_1)}(X,M_{^2{\bf k}^2})\|\\
		\end{array}\right\}.
	$$
	Крім того, $M_{{\bf k}} \in {\rm int}A_{{\bf k}}(X)$ тоді і тільки тоді, коли $M_{{\bf k}^2} \in {\rm int} A_{\K^2}(X)$ і 
	$M_{k_1}> \|\phi^{(k_1)}(X,M_{{\bf k}^2})\|$ (при непарному $d$) або 	$M_{k_1}> \|\phi^{(k_1)}(X,M_{^2{\bf k}^2})\|$ (при парному $d$).
\end{theorem}
{\bf Доведення.} Достатньо довести випадок, коли $X=L_{\infty,\infty}^{r,\cup}(\RR_-)$. Для скорочення записів будемо опускати позначення класу $X$. 

Доведемо необхідність вказаних умов. Нехай 
\begin{equation}\label{MkisAdmissible}
	M_{{\bf k}} \in A_{{\bf k}}. 
\end{equation}
 Оскільки множина $A_{\K^2}$ є опуклою, то $A_{\K^2} = {\rm int}A_{\K^2} \bigcup \left(\partial A_{\K^2}\bigcap A_{\K^2}\right)$. Це означає, що $M_{\K^2}\in {\rm int}A_{\K^2}$ або $M_{\K^2}\in \partial A_{\K^2}\bigcap A_{\K^2}$. Нехай  
\begin{equation}\label{**}
	M_{{\bf k}^2} \in {\rm int} A_{\K^2}.
\end{equation}
Необхідність умови
\begin{equation}\label{***}
	M_{k_1}\ge \|\phi^{(k_1)}(M_{{\bf k}^2})\|
\end{equation}
у випадку непарного $d$ і умови  
\begin{equation}\label{even}
	M_{k_1} > \|\phi^{(k_1)}(M_{^2{\bf k}^2})\|
\end{equation}
у випадку парного $d$ випливає з леми~\ref{l::1} --- у випадку парного $d$ слід ще зауважити, що з~\eqref{**} та леми~\ref{l::1} слідує, що виконується строга нерівність
$$\|\phi^{(r)}(M_{^2{\bf k}^2})\| < M_{k_d}.$$

Нехай 
\begin{equation}\label{Mk2isOnBoundary}
	M_{{\bf k}^2} \in \partial A_{\K^2}(X)\bigcap A_{\K^2}(X).
\end{equation}
У випадку непарного $d = 2n+1$ сплайн $\phi(M_{\K^2})$ має не більше ніж $n-1$ вузол, у випадку парного $d = 2n$ сплайн $\phi(M_{^2\K^2})$ --- не більше ніж $n-2$ вузла. Нехай функція $x\in L^{r,\cup}_{\infty,\infty}({\RR_-})$ така, що $M_{\K^2}(x)=M_{\K^2}$ (у випадку непарного $d$) або $M_{^2\K^2}(x)=M_{^2\K^2}$ (у випадку парного $d$). Тоді в силу леми~\ref{l::1} $x^{(k_2)}\equiv\phi^{(k_2)}(M_{\K^2})$ у випадку непарного $d$ та  $x^{(k_2)}\equiv\phi^{(k_2)}(M_{^2\K^2})$  у випадку парного $d$.  

У випадку непарного $d$, коли $k_1>0$, це означає, що  $x^{(k_1)}\equiv\phi^{(k_1)}(M_{\K^2})$  (а отже $\left\|x^{(k_1)}\right\|=\left\| \phi^{(k_1)}(M_{{\bf k}^2})\right\|$). Якщо ж $k_1=0$, то $x\equiv\phi(M_{\K^2})+C$, $C\geq 0$ (а отже $\left\|x\right\|\geq\left\| \phi(M_{{\bf k}^2})\right\|$). Аналогічно у випадку парного $d$. 

У випадку, коли $d$ непарне, має місце~\eqref{**} і у нерівності~\eqref{***} має місце рівність, то 
\begin{equation}\label{MkisOnBoundary}
	M_{\K}\in \partial A_{\K}\bigcap A_{\K},
\end{equation}
оскільки з леми~\ref{l::1} слідує, що для всіх $\varepsilon > 0$ $(M_{k_1}-\varepsilon,M_{k_2},\dots,M_{k_{d}})\notin A_{\K}$. Аналогічно у випадку, коли $d$ парне, має місце~\eqref{**} і в нерівності~\eqref{even} має місце рівність, то справедливе включення~\eqref{MkisOnBoundary}.

Таким чином, необхідність наведених умов доведено.

Доведемо тепер достатність вказаних умов.

Достатність вказаних умов у випадку, коли виконується~\eqref{Mk2isOnBoundary} та $d$ непарне очевидна. У випадку парного $d$ відмітимо, що з~\eqref{Mk2isOnBoundary} та леми~\ref{l::1} слідує, що $\|\phi^{(r)}(M_{^2{\bf k}^2})\| = M_{k_d}$.

Доведемо достатність у випадку виконання умови~\eqref{**}. 

Розглянемо випадок непарного $d=2n+1$. Достатність очевидна у випадку $k_1 = 0$, тому далі будемо вважати, що $k_1 > 0$.

Можемо вважати, що в~\eqref{***} має місце строга нерівність. 
З~\eqref{**} слідує, що для всіх $\varepsilon > 0$ $(\|\phi(M_{\K^2})\|+\varepsilon,M_{k_2},\dots, M_{k_d})\in {\rm int}A_{{\bf K}}$, де ${\bf K} = (0,k_2,k_3,\dots, k_d)$. Нехай $a_n$ --- самий правий вузол сплайна $\phi(M_{\K^2})$. В силу леми~\ref{l::2} для всіх $\varepsilon > 0$ існує $L_{\infty,\infty}^{r,\cup}(\RR_-)$ -- ідеальний сплайн $\psi = \psi(\varepsilon)=\phi({\bf a}(\varepsilon),{\bf \lambda}(\varepsilon))$, де вектори ${\bf a}(\varepsilon)=(a_1(\varepsilon),a_2(\varepsilon),\dots,a_{n+1}(\varepsilon))$ і ${\bf \lambda(\varepsilon) }= (\lambda_1(\varepsilon),\dots,\lambda_{n+1}(\varepsilon))$ такі, що  $a_1(\varepsilon)>a_2(\varepsilon)>\ldots>a_{n+1}(\varepsilon)= \frac {a_n}{2}$ і $\lambda_1(\varepsilon),\dots,\lambda_{n+1}(\varepsilon)>0$, такий, що для $i=2,\dots, d$ $\left\|\psi^{(k_i)}(\varepsilon)\right\| = M_{k_i}$ і $\left\|\psi(\varepsilon)\right\| = \varepsilon + \left\|\phi(M_{\K^2})\right\|$. Для $i=1,\dots,n+1$ покладемо $a_i^*:=\lim\limits_{\varepsilon\to+\infty}a_i(\varepsilon)$ (при цьому границі можуть бути скінченними або нескінченними). Ясно, що $a_1^*=\infty$. Припустимо, що функція $\lambda_1(\varepsilon)a_1^{r-k_1}(\varepsilon)$ обмежена при  $\varepsilon \to+\infty$. Тоді в силу леми~\ref{l::-1} $\lambda_1(\varepsilon)a_1^{r-k_i}(\varepsilon)\to 0$ при $\varepsilon \to+\infty$ для всіх $i=2,\dots,d$. Це означає, що $a_2^*=\infty$, оскільки в силу зауваження~\ref{z::3} і того факту, що $a_{n+1}(\varepsilon)\equiv \frac {a_n}{2}$, всі границі $a_i^*$, $i=2,\dots, n+1$, не можуть бути скінченними. Повторюючи аналогічні міркування ми отримаємо, що для деякого $i=1,\dots,n$ функція $\lambda_i(\varepsilon)a_i^{r-k_1}(\varepsilon)$ є необмеженою при $\varepsilon\to\infty$. Це означає, що норму $\left\|\psi^{(k_1)}(\varepsilon)\right\|$ можна зробити скільки завгодно великою. В силу опуклості множини $A_{{\bf k}}$ це означає, що виконується~\eqref{MkisAdmissible}.

Крім того, оскільки при виконанні~\eqref{**} для довільного $M_{k_1}$, для якого виконується строга нерівність~\eqref{***}, маємо~\eqref{MkisAdmissible}, то для таких $M_{k_1}$ $M_{\K}\in {\rm int}A_{{\bf k}}$. 

Нехай тепер $d = 2n$. 
В силу леми~\ref{l::2} для всіх $a>0$ існує сплайн $\psi(a)$ з $n$ вузлами, один з яких розташовану у точці  $a$ такий, що виконуються рівності $M_{\K^2}(\psi(a)) = M_{\K^2}$. При $a\to 0$ має місце поточкова збіжність $\psi(a)\to\phi(M_{^2{\bf k}^2})$, а отже при достатньо малих $a$ $\left\|\psi^{(k_1)}(a)\right\|<M_{k_1}$. Переходячи до границі при $a\to\infty$ і застосовуючи міркування аналогічні до випадку непарного $d$ отримаємо справедливість~\eqref{MkisAdmissible}.

Теорему доведено.

\begin{remark}
	Теорема~\ref{main} є аналогом результатів $\S 6$ глави 2 в~\cite{Karlin66}.
\end{remark}
\bibliographystyle{my_ugost2003s}
\bibliography{bibliography}

\end{document}